\documentclass[11pt, a4paper]{amsart}
\usepackage[top=35truemm,bottom=35truemm,left=35truemm,right=35truemm]{geometry}
\usepackage{setspace} 
\usepackage{amsmath, amsfonts,amsthm,amssymb,mathrsfs,amscd}
\usepackage{enumerate}
\usepackage[dvipdfmx]{graphicx}
\usepackage{ascmac}
\usepackage[arrow,matrix]{xy}
\usepackage{color}
\usepackage{array}
\usepackage{pifont}
\usepackage{longtable}
\usepackage{here}
\usepackage{braket}
\usepackage{mathtools}

\theoremstyle{definition}
\newtheorem{theo}{Theorem}
\newtheorem*{theo*}{Theorem}
\newtheorem{defi}{Definition}
\newtheorem{exa}{Example}
\newtheorem{prop}{Proposition}
\newtheorem{cor}{Corollary}
\newtheorem{lem}{Lemma}

\def\C{\mathbb{C}}

\def\N{\mathbb{N}}
\def\H{\mathbb{H}}
\def\R{\mathbb{R}}

\begin{document}

\title[ ]{Construction of infinite series of non-simple ideal hyperbolic Coxeter 4-polytopes and their growth rates}
\author{Tomoshige Yukita}
\address{Department of Mathematics, School of Education, Waseda University, Nishi-Waseda 1-6-1, Shinjuku, Tokyo 169-8050, Japan}
\email{yshigetomo@suou.waseda.jp}
\subjclass[2010]{Primary~20F55, Secondary~20F65}
\keywords{Coxeter group; growth function; growth rate; Perron number}
\date{}
\thanks{}

\begin{abstract}
In this paper, we construct infinite series of non-simple ideal hyperbolic Coxeter 4-polytopes whose growth rates are Perron numbers. 
This infinite series is the first example of such a non-compact infinite polytopal series.
\end{abstract}

\maketitle


\setstretch{1.1}
\section{Introduction}
Let $\mathbb{H}^{d}$ denote the upper half-space model of hyperbolic $d$-space and $\overline{\mathbb{H}}^{d}$ its closure in $\mathbb{R}^d\cup{\{\infty\}}$. 
A convex polytope $P\subset{\overline{\mathbb{H}}^{d}}$ of finite volume is called a \textit{Coxeter polytope} if all of its dihedral angles are of the form $\frac{\pi}{k}$ for an integer $k\geq{2}$ or $k=\infty$, i.e. the intersection of respective facets is a point on the boundary $\partial{\mathbb{H}^{d}}$.
The set $S$ of reflections with respect to facets of $P$ generates a discrete group $\Gamma$, called a \textit{hyperbolic Coxeter group},
and the pair $(\Gamma, S)$ is called the \textit{Coxeter system} associated with $P$. 
Then $P$ becomes a fundamental domain for $\Gamma$. 
If $P$ is compact (resp. non-compact), the hyperbolic Coxeter group $\Gamma$ is called \textit{cocompact} (resp. \textit{cofinite}). 
The \textit{growth series} $f_{S}(t)$ of $(\Gamma, S)$ is the formal power series $\sum_{l=0}^{\infty}a_{l}t^{l}$ where $a_{l}$ is the number of elements of $\Gamma$ whose word length with respect to $S$ is equal to $l$. 
Then $\tau_{\Gamma} :=\limsup_{l \rightarrow \infty} \sqrt[l]{a_l}$ is called the \textit{growth rate} of $(\Gamma, S)$. 
By means of the Cauchy-Hadamard theorem, $\tau_{\Gamma}$ is equal to the reciprocal of the radius of convergence $R$ of $f_{S}(t)$. 
The growth series and the growth rate of a hyperbolic Coxeter polytope $P$ is defined to be the growth series and the growth rate of the Coxeter system $(\Gamma, S)$ associated with $P$, respectively. 
It is known that the growth rate of a hyperbolic Coxeter polytope is a real algebraic integer bigger than 1 \cite{dlH1}. 
Recall that a real algebraic number $\tau>1$ is a Perron number, if and only if all of its other algebraic conjugates are less than $\tau$ in absolute value.
It is known that the growth rates of 2 and 3-dimensional hyperbolic Coxeter polytopes are always Perron numbers (\cite{CW}, \cite{F}, \cite{Pa}, \cite{Y1}, \cite{Y2}). 
From now on, we consider the growth rates of hyperbolic Coxeter 4-polytopes. 
In the study of the growth rates of compact hyperbolic Coxeter 4-polytopes, (1) Kellerhals and Perren showed that the growth rates of compact hyperbolic Coxeter 4-polytopes with at most 6 facets are Perron numbers \cite{KP} and (2) T.Zehrt and C.Zehrt \cite{Z} and Umemoto \cite{U} constructed infinite series of compact hyperbolic Coxeter 4-polytopes and proved that their growth rates are 2-Salem numbers which are particular Perron numbers. 
In this paper, we consider new infinite series of {\it ideal} and non-simple hyperbolic Coxeter 4-polytopes and prove that their growth rates are Perron numbers. 
In this way, we provide the first example of such a non-compact infinite polytopal series and prove that their growth rates are Perron numbers.

The organization of the present paper is as follows.
In Section 2, we review useful formulas which allow us to calculate the growth function of a hyperbolic Coxeter polytope. 
In Section 3, we explain a method to determine the distribution of roots of a real polynomial. 
Then, we construct infinite series $\{P_n\}_{n\in{\N}}$ of non-simple ideal hyperbolic Coxeter 4-polytopes in Section 4. 
Finally, we apply the method introduced in Section 3 to the denominator polynomial $D_n(t)$ of the growth functions $f_n(t)$ of the polytope $P_n$ in Section 5. 
In Appendix, we list numerical data of the denominator polynomials of the growth functions.

\section{preliminaries}
In this Section,  we introduce the relevant notation and review Solomon's and Steinberg's formulas in order to calculate the growth functions of hyperbolic Coxeter polytopes. 

\begin{defi}{\rm(Coxeter system, Coxeter diagram, growth rate)}

(i) A \textit{Coxeter system} $(\Gamma,S)$ consists of a group $\Gamma$ and a finite set of generators $S\subset{\Gamma}$, $S=\{s_i\}_{i=1}^N$, with relations $(s_is_j)^{m_{ij}}$ for each $i , j$ , where $m_{ii}=1$ and $m_{ij}\geq{2}$ or $m_{ij}=\infty$ for $i\neq{j}$.
We call $\Gamma$ a \textit{Coxeter group}.
For any subset $I\subset{S}$, we define $\Gamma_I$ to be the subgroup of $\Gamma$ generated by $\{s_i\}_{i\in{I}}$. 
Then $(\Gamma_I, I)$ is a Coxeter system in its own right and $\Gamma_I$ is called the \textit{Coxeter subgroup} of $\Gamma$ generated by $I$.

(ii) The \textit{Coxeter diagram} $X(\Gamma, S)$ of $(\Gamma,S)$ is constructed as follows: \\
Its vertex set is $S$. 
If $m_{ij}\geq{3} \ (s_i\neq{s_j} \in{S})$, we join the pair of vertices by an edge. 
For each edge, we label it with $m_{ij}$ if $m_{ij}\geq{4}$. 
Note that the Coxeter diagram of $(\Gamma_I, I)$ for each subset $I\subset{S}$ is a subdiagram of $X(\Gamma, S)$. 

(iii) The \textit{growth series} $f_{S}(t)$ of $(\Gamma,S)$ is the formal power series $\sum_{l=0}^{\infty}a_lt^l$ where $a_l$ is the number of elements of $\Gamma$ whose word length with respect to $S$ is equal to $l$. Then $\tau_{(\Gamma, S)} =\limsup_{l \rightarrow \infty} \sqrt[l]{a_l}$ is called the \textit{growth rate} of $(\Gamma,S)$.

\end{defi}

A Coxeter system $(\Gamma, S)$ is \textit{irreducible} if the Coxeter diagram of $(\Gamma,S)$ is connected. 
We recall Solomon's formula and Steinberg's formula which enable us to express the growth series of Coxeter systems as rational functions. 

\begin{theo}\label{Solomon}{\rm (Solomon's formula)}\cite{So} 
\textit{The growth series $f_S(t)$ of an irreducible finite Coxeter system $(\Gamma, S)$ can be written as 
$f_S(t)=[m_1+1, m_2+1, \cdots, m_p+1]$ where $[n] =1+t+ \cdots +t^{n-1}, [m, n]=[m][n]$,etc., and where $\{m_1, m_2, \cdots, m_p \}$
is the set of exponents of $(\Gamma, S)$.}
\end{theo}

The exponents of irreducible finite Coxeter groups are shown in Table 1 (see \cite{Hu} for details).

\begin{table}[h]
\begin{center}
\caption{Exponents}
\begin{tabular}{|c|c|c|}
\hline 
Coxeter group & Exponents & growth series \\ 
\hline
$A_n$ & $1,2,\cdots,n$ & $[2,3,\cdots,n+1]$ \\ 
\hline
$B_n$ & $1,3,\cdots,2n-1$ & $[2,4,\cdots,2n]$ \\ 
\hline
$D_n$ & $1,3,\cdots,2n-3,n-1$ & $[2,4,\cdots,2n-2][n]$ \\ 
\hline
$E_6$ & 1,4,5,7,8,11 & [2,5,6,8,9,12] \\ 
\hline
$E_7$ & 1,5,7,9,11,13,17 & [2,6,8,10,12,14,18] \\ 
\hline
$E_8$ & 1,7,11,13,17,19,23,29 & [2,8,12,14,18,20,24,30] \\ 
\hline
$F_4$ & 1,5,7,11 & [2,6,8,12] \\ 
\hline
$H_3$ & 1,5,9 & [2,6,10] \\ 
\hline
$H_4$ & 1,11,19,29 & [2,12,20,30] \\ 
\hline
$I_2(m)$ & 1,$m-1$ & [2,$m$] \\ 
\hline
\end{tabular}
\end{center}
\end{table}

\begin{theo}\label{Steinberg}{\rm (Steinberg's formula)}\cite{St} 
\textit{Let $(\Gamma, S)$ be an infinite Coxeter system.  
Set $\mathcal{F}:=\Set{I\subset{S}|\Gamma_I \textrm{ is a finite Coxeter subgroup of }\Gamma}$. 
Denote by $f_I(t)$ the growth series of the Coxeter system $(\Gamma_I, I)$ for each $I\subset{S}$.
Then
$$
\frac{1}{f_S(t^{-1})}=\sum _{I \in \mathcal{F}} \frac{(-1)^{|I|}}{f_I(t)}.
$$}
\end{theo}

By Theorem \ref{Solomon} and Theorem \ref{Steinberg}, the growth series of $(\Gamma, S)$ is represented by a rational function $\frac{p(t)}{q(t)} (p,q\in{\mathbb{Z}[t]})$. 
The rational function $\frac{p(t)}{q(t)}$ is called the \textit{growth function} of $(\Gamma, S)$. 
The radius of convergence $R$ of the growth series $f_{S}(t)$ is equal to the positive real root of $q(t)$ which has the smallest absolute value among all the roots of $q(t)$.

In this paper, we are interested in Coxeter groups which act discontinuously on hyperbolic space $\H^d$.

\begin{defi}(Upper half-space model of hyperbolic $d$-space)

The upper half-space $\H^d=\Set{(x_1, \cdots, x_d)\in{\R^d} | x_d>0}$ equipped with the metric $\frac{|dx|}{x_d}$ is a model of hyperbolic $d$-space, so called the \textit{upper half-space model}. 
The boundary $\partial{\H^d}$ of $\H^d$ in the one-point compactification $\R^d\cup{\{\infty\}}$ of Euclidean $d$-space $\R^d$ is called the \textit{boundary at infinity}. 
We denote the closure of a subset $A\subset{\R^d\cup{\{\infty\}}}$ by $\bar{A}$. 
\end{defi}

By identifying $\R^{d-1}$ with $\R^{d-1}\times{\{0\}}$ in $\R^d$, 
the boundary at infinity $\partial{\H^d}$ is equal to $\R^{d-1}\cup{\{\infty\}}$. 
A subset $H\subset{\H^d}$ is called a \textit{hyperplane} of $\H^{d}$ 
if and only if it is a Euclidean hemisphere or a half-plane orthogonal to $\R^{d-1}$. 

\begin{defi}{\rm(hyperbolic polytope)} 

A subset $P\subset{\overline{\mathbb{H}}^d}$ is called a \textit{hyperbolic polytope} 
if $P$ can be written as the intersection of finitely many closed half-spaces: 
$P=\cap{H^-_i}$, where $H^-_i$ is the closed domain of $\mathbb{H}^d$ bounded by a hyperplane $H_i$. 
\end{defi}

Suppose that $H_i\cap{H_j}\neq{\emptyset}$  in $\mathbb{H}^{d}$. 
Then we define the {\textit {dihedral angle}} between $H_i$ and $H_j$ as follows: let us choose a point $x\in{H_i\cap{H_j}}$ and consider the outer normal vectors $u_i$ and $u_j$. 
Then the dihedral angle between $H_i$ and $H_j$ is defined as the real number $\theta\in{[0, \pi)}$ satisfying $\cos{\theta}=-(u_i, u_j)$ where $(\cdot , \cdot)$ denotes the Euclidean inner product on $\mathbb{R}^d$ at $x$.

If $\overline{H_i}\cap{\overline{H_j}}\in{\overline{\mathbb{H}}^{d}}$ is a point on $\partial{\mathbb{H}^d}$, 
then we define the dihedral angle between $H_i$ and $H_j$ to be equal to zero.

\begin{defi}{\rm(hyperbolic Coxeter polytope)}

A hyperbolic polytope $P\subset{\overline{\mathbb{H}}^{d}}$ of finite volume is called a {\textit{hyperbolic Coxeter polytope}} 
if all of its dihedral angles have the form $\frac{\pi}{k}$ for an integer $k\geq{2}$ or $k=\infty$ 
if the intersection of respective bounding hyperplanes is a point on $\partial{\mathbb{H}^{d}}$.
\end{defi}

Notice that a hyperbolic polytope in $\overline{\mathbb{H}}^{d}$ is of finite volume 
if and only if it is the convex hull of finitely many points in $\overline{\H}^d$. 
If $P\subset{\overline{\mathbb{H}}^{d}}$ is a hyperbolic Coxeter polytope, 
the set $S$ of all reflections with respect to facets of $P$ generates a discrete group $\Gamma$. 
It is known that $(\Gamma, S)$ is a Coxeter system, so that $\Gamma$ is a Coxeter group. 
We call $\Gamma$ the \textit{$d$-dimensional hyperbolic Coxeter group}, 
and the pair $(\Gamma, S)$ is called the \textit{Coxeter system associated with $P$}. 
In the sequel, the growth function and the growth rate of the Coxeter system $(\Gamma, S)$ associated with $P$ 
are called the \textit{growth function of $P$} and the \textit{the growth rate of $P$}. 
The growth function and the growth rate of $P$ are denoted by $f_P(t)$ and $\tau_{P}$.

\begin{defi}(Gram matrix, Coxeter scheme)
Let $P=\cap_{i=1}^{N}H_i^{-}$ be a hyperbolic Coxeter polytope. 
To every pair of hyperplanes $H_i$ and $H_j$, 
define 
\[
c_{ij}=\begin{cases}
        1 & \text{if }i=j \\
        -\cos{\frac{\pi}{m_{ij}}} & \text{if they intersect at the dihedral angle }\frac{\pi}{m_{ij}} \\
        -1 & \text{if its intersection is a point on } \partial{\H}^d \\
        -\cosh{d(H_i, H_j)} & \text{if they do not intersect}
        \end{cases}
\]
where, $d(H_i, H_j)$ is the hyperbolic distance between them. 
The $N\times{N}$ symmetric matrix $M(P)=(c_{ij})$ is called the \textit{Gram matrix of $P$}. 
The \textit{Coxeter scheme $X(P)$ of $P$} is defined as follows; 
Its vertex set is $\{H_1, \cdots, H_N\}$. 
If the dihedral angle $\frac{\pi}{m_{ij}}$ between hyperplanes $H_i$ and $H_j$ is less than $\frac{\pi}{2}$, 
we join the pair of vertices by an edge. 
For each edge, we label it with $m_{ij}$ if $m_{ij}\geq{4}$. 
Two vertices are joined by a dotted edge labeled with the hyperbolic distance between corresponding hyperplanes if they do not intersect.
\end{defi}
A subscheme of a Coxeter scheme $X(P)$ is called \textit{elliptic} (resp. \textit{parabolic}) 
if the corresponding submatrix of the Gram matrix $M(P)$ is positive definite (resp. positive semi-definite and its rank equals $d-1$). 
Note that elliptic subschemes correspond to finite Coxeter systems. 

\begin{theo}(Theorem 2.2, p.109 and Theorem 2.5, p.110 \cite{V})
Given a hyperbolic Coxeter polytope $P$, the faces (resp. vertices of infinity) of $P$ correspond to the elliptic (resp. parabolic) subschemes of the Coxeter scheme of $P$. 
\end{theo}

\section{Method for deciding the distribution of the roots of a real polynomial}
In this Section, we review Sturm's theorem and Kronecker's theorem.
Sturm's theorem shows how we can determine the distribution of real roots of a real polynomial 
and Kronecker's theorem tells us how to count roots of a real polynomial 
contained in a closed disk of radius $r$ centered at the origin $0$ in the complex plane $\C$. 
The argument in the present section is based on \cite{C}, \cite{K} and \cite{P}.
\subsection{Sturm's theorem}
\begin{defi}(Sturm sequence)
Let $f$ and $g$ be real polynomials. 
We may assume that $\deg{f}\geq{\deg{g}}$. 
By the Euclidean algorithm, we define polynomials $f_2, \cdots, f_r$ as follows: 
\begin{eqnarray*}
f &=& q_1g-f_2, \hspace{5mm} \deg{g}>\deg{f_2}. \\
g &=& q_2f_2-f_3, \hspace{5mm} \deg{f_2}>\deg{f_3}. \\
f_2 &=& q_3f_3-f_4, \hspace{5mm} \deg{f_3}>\deg{f_4}. \\
\vdots && \\
f_{r-2} &=& q_{r-1}f_{r-1}-f_r, \hspace{5mm} \deg{f_{r-1}}>\deg{f_{r}}. \\
f_{r-1} &=& q_rf_r. 
\end{eqnarray*} 
Then, the finite sequence $f_0:=f, f_1:=g, f_2, \cdots, f_r$ of real polynomials is called the \textit{Sturm sequence of $f$ and $g$}.
\end{defi}

Note that $f_r$ is the greatest common divisor of polynomials $f$ and $g$. 
For any $t_0\in{\R}$, the number of sign changes in the Sturm sequence of $f$ and $g$ at $t_0$ is denoted by $w(t_0)$, 
that is, 
$w(t_0)$ is the number of sign changes in the sequence $f(t_0), g(t_0), f_2(t_0), \cdots, f_r(t_0)$ ignoring zeros.

\begin{exa}
Let $f(z):=z^5-3z-1$ and $g(z):=f'(z)=5z^4-3$. 
Then, the Sturm sequence of $f$ and $g$ can be calculated as follows: 
\begin{eqnarray*}
f(z)&=&z^5-3z-1. \\
g(z)&=&5z^4-3. \\
f_2(z)&=&12z+5. \\
f_3(z)&=&1.
\end{eqnarray*}
We consider the number of sign changes in the Sturm sequence at $-2$. 
We have $f(-2)=-27, g(-2)=77, f_2(-2)=-19, f_3(-2)=1$, so that $w(-2)$ is equal to 3. 
\end{exa}

\begin{theo} (Sturm's theorem, Theorem8.8.15 \cite{C})
\textit{Let $f$ be a real polynomial and $f_0, f_1, \cdots, f_r$ be the Sturm sequence of $f$ and $f'$. 
Suppose that $a, b\in{\R}$ are not roots of $f$ and $a<b$. 
Then the number of distinct real roots of $f$ in the closed interval $[a, b]$ is equal to $w(a)-w(b)$.} 
\end{theo}

From now on, we assume that real polynomials $f$ and $g$ have no common roots. 
For each real root $t_0$ of $f$, the number of sign changes in $f, g$ satisfies one of the following three conditions; 
\begin{itemize}
\item[(i)] the number of sign changes in $f(t), g(t)$ decreases by 1 when $t$ pass through $t_0$. 
\item[(ii)] the number of sign changes in $f(t), g(t)$ increases by 1 when $t$ pass through $t_0$.
\item[(iii)] the number of sign changes in $f(t), g(t)$ does not vary when $t$ pass through $t_0$.
\end{itemize}
We assign the number $\varepsilon_{t_0}=1, -1$ and $0$ to each root $t_0$ of $f$ 
when the number of sign changes of $f$ and $g$ satisfies the condition (i), (ii) and (iii), respectively.  
The following theorem is proved analogously to Sturm's theorem.

\begin{theo} 
\textit{Suppose that real numbers $a$ and $b$ are not roots of $f$. 
Then, the following identity holds for the Sturm sequence of $f$ and $g$. 
\[
\sum_{t_0\in{[a, b]}:\text{ a zero of }f}\varepsilon_{t_0}=w(a)-w(b).
\]}
\end{theo}

\subsection{Separation of complex roots} 
We use the following notation in this Subsection: 
\begin{itemize}
\item $\C_z$ and $\C_w$ denote respectively complex plane with coordinate $z=x+iy$ and $w=u+iv$. 
\item $S_r\subset{\C_z}$ is a circle of radius $r>0$ centered at the origin $0\in{\C_z}$. 
\item $B_r\subset{\C_z}$ is a open disk of radius $r>0$ centered at $0$. 
\item  A parameter for $S_r$ is given as follows: 
\[
z(t)=r\dfrac{t^2-1}{t^2+1}-ir\dfrac{2t}{t^2+1}.
\]
\item $f(z)$ is a real polynomial of a complex variable $z$. 
\item By using two real polynomials $\varphi_r(t), \psi_r(t)$ of a real variable $t$, $f(z(t))$ can be represented as
\[
f(z(t))=\dfrac{\varphi_r(t)+i\psi_r(t)}{(t^2+1)^{\deg{f}}} \text{ on }S_r. 
\]
\end{itemize}

\begin{lem}
\textit{Suppose that $f(z)$ has no roots on $S_r$. 
Given $M>0$ s.t. the closed interval $[-M, M]$ contains all real roots of $\varphi_r$, the following identity holds for the Sturm sequence of $\varphi_r$ and $\psi_r$. 
\[
\sum_{t_0\in{[-M, M]}:\text{ a zero of }\varphi_r}\varepsilon_{t_0}=w(-M)-w(M).
\]}
\end{lem}
\proof The assumption that $f(z)$ has no roots on $S_r$ implies that 
real polynomials $\varphi_r(t)$ and $\psi_r(t)$ do not have common real roots. 
Therefore, we can apply Theorem 4 to $\varphi_r(t)$ and $\psi_r(t)$. \qed

\begin{figure}[htbp]
\begin{center}
 \includegraphics [scale=.5, clip]{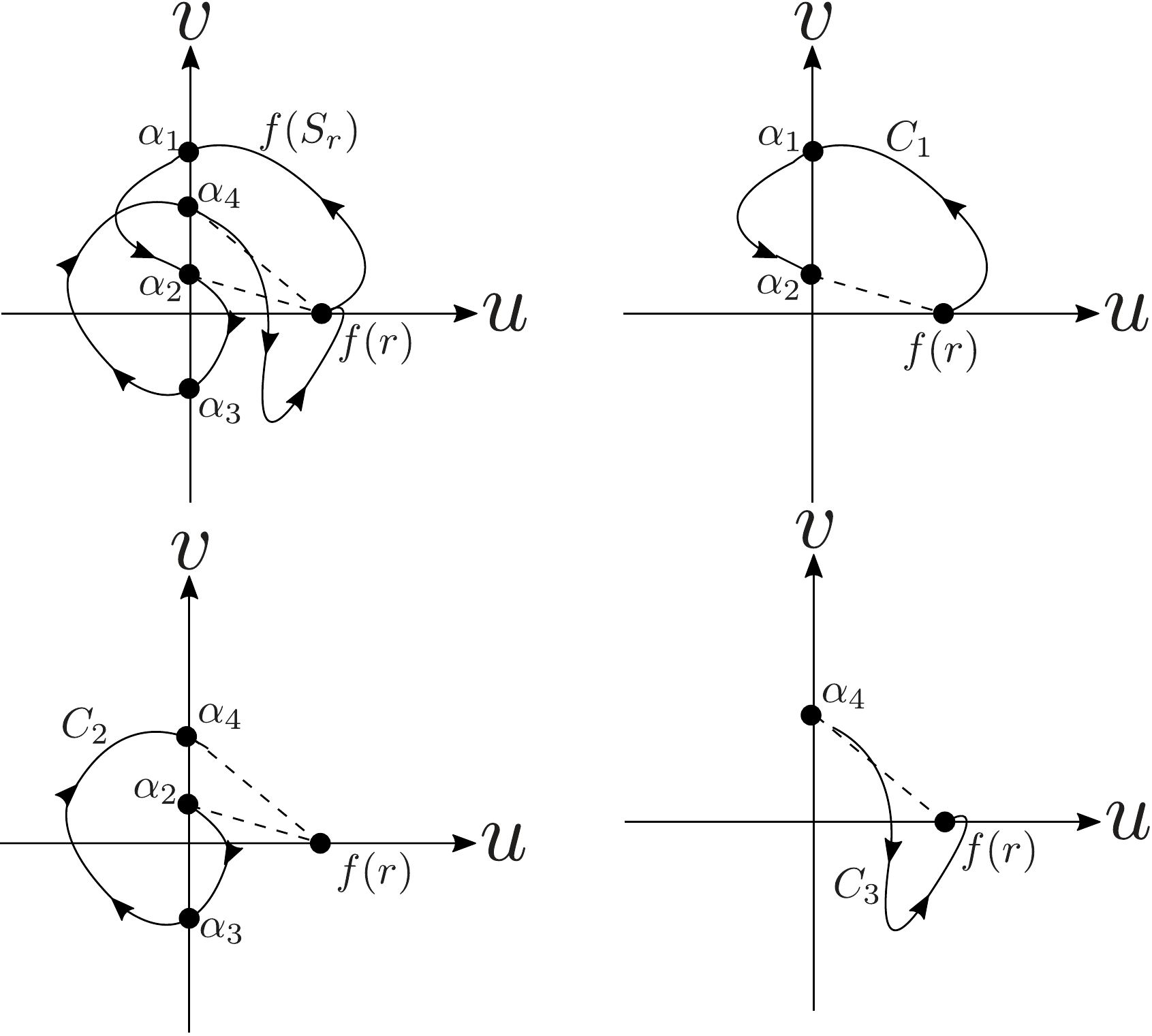}
\end{center}
\caption{Division of the closed curve $f(S_r)$}
\label{figure1}
\end{figure}

\begin{figure}[htbp]
\begin{center}
 \includegraphics [scale=.5, clip]{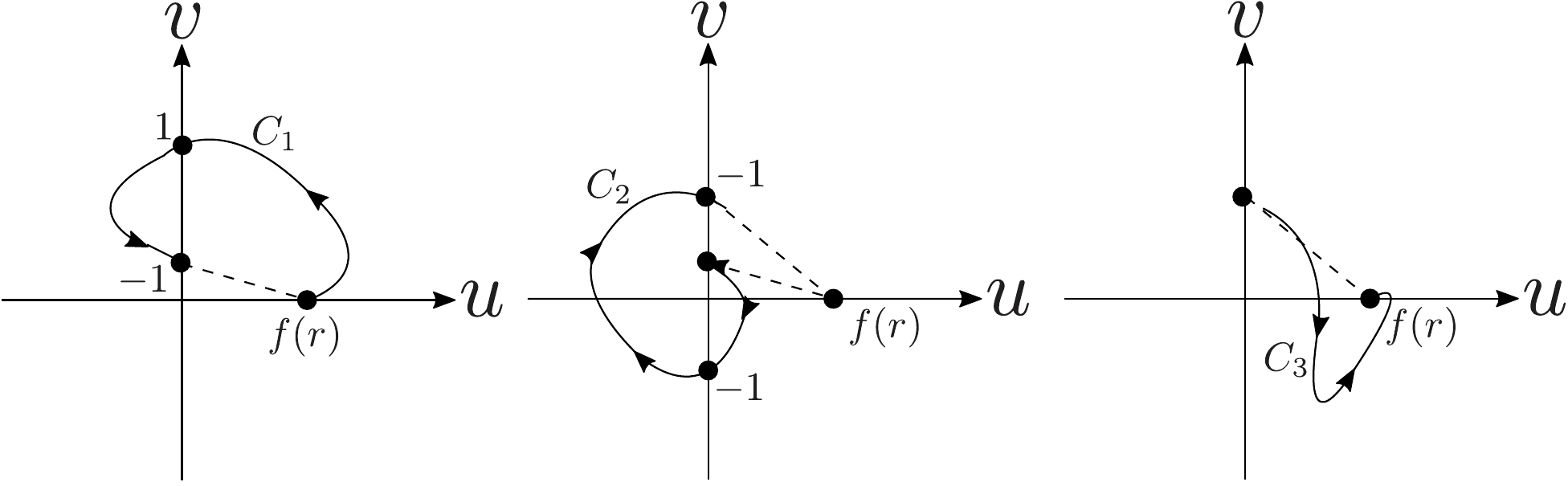}
\end{center}
\caption{Assigning the number to crossing points}
\label{figure2}
\end{figure}

By identifying $f(z)$ as a holomorphic function from $\C_z$ to $\C_w$, 
we give a parameter for the closed curve $f(S_r)$ as 
$w(t)=\frac{\varphi_r(t)}{(t^2+1)^{\deg{f}}}+i\frac{\psi_r(t)}{(t^2+1)^{\deg{f}}}$. 
In order to calculate the winding number of $f(S_r)$, 
we divide $f(S_r)$ into closed curves $C_1, \cdots, C_m$ as follows; 
trace $f(S_r)$ from the initial point $f(r)=\lim_{t\to{-\infty}}w(t)$, 
and if the curve crosses the $v$-axis twice, then we mark each crossing point with $\alpha_1$ and $\alpha_2$ 
and go back to the initial point $f(r)$ along the straight line from the point $\alpha_2$ to the initial point $f(r)$. 
This locus makes the closed curve $C_1$. 
After that, we go back to $f(S_r)$ along the straight line from $f(r)$ to $\alpha_2$. 
By repeating this procedure, the closed curve $f(S_r)$ is divided into closed curves $C_1, \cdots, C_m$ (see Fig \ref{figure1}).
Under the division of $f(S_r)$, the winding number of $f(S_r)$ equals to the sum of the winding numbers of closed curves $C_1, \cdots, C_m$. 
To calculate the winding number of each closed curve $C_i$, 
we assign the number $\chi_{\alpha_k}=1$ (resp. $\chi_{\alpha_k}=-1$) to a crossing point $\alpha_{k}$ of the $v$-axis and $C_i$ 
if the argument of $C_i$ is increasing (resp. decreasing) around the crossing point $\alpha_k$ (see Fig \ref{figure2}). 
Then, the winding number of $C_i$ is equal to the sum of $\frac{1}{2}\chi_{\alpha_k}$ on each crossing point $\alpha_k$.
Note that if $C_i$ has no crossing points of the $v$-axis and $C_i$, then the winding number of $C_i$ is equal to $0$. 
For example, the winding number of $C_1, C_2$ and $C_3$ in Fig \ref{figure2} is equal to $0, -1$ and $0$, respectively.  
This observation shows that the winding number of $f(S_r)$ is equal to the sum of the number $\frac{1}{2}\chi_{\alpha}$ on each crossing point $\alpha$ of the $v$-axis and $f(S_r)$.

Let us now consider the Sturm sequence of polynomials $\varphi_r(t)$ and $\psi_r(t)$. 
Every crossing point of the curve $f(S_r)$ corresponds to a root of $\varphi_r(t)$. 
For any root $t_0\in{\R}$ of $\varphi_r(t)$, the argument of $f(S_r)$ is increasing (resp. decreasing) if $\varepsilon_{t_0}=-1$ (resp. $\varepsilon_{t_0}=1$). 
This observation, together with Theorem 4 and the argument principle, implies the following equalities. 
\begin{eqnarray*}
\#\Set{z\in{B_r} | z\text{ is a root of }f(z)} &=& \text{The winding number of }f(S_r) \\
                                               &=& \frac{1}{2}\sum_{\alpha_{k}:\text{ a mark on }F(S_r)} \chi_{\alpha_k} \\
                                               &=& \frac{1}{2}\sum_{t_0:\text{ a zero of }\varphi_r}-\varepsilon_{t_0}. 
\end{eqnarray*}
By Lemma 1, we obtain Kronecker's theorem.

\begin{theo} (Kronecker's theorem, Theorem1.4.6 \cite{P})
\textit{Suppose that $f(z)$ have no roots on $S_r$. 
Then the number of roots of $f$ contained in $B_r$ equals to $\frac{w(M)-w(-M)}{2}$, 
where $M>0$ is a real number such that $[-M, M]$ contains all roots of $\varphi_r(t)$. }
\end{theo}

If we substitute $z(t)=r\dfrac{t-i}{t+i}$ for $f(z)$, then $f(z(t))$ can be rewritten as follows:
\[
f(z(t))=\dfrac{\Phi(t)+i\Psi(t)}{(t+i)^{\deg{f}}}.
\]
Since $\arg{f(z(t))}=\arg{\{\Phi(t)+i\Psi(t)\}}-\arg{(t+i)^{\deg{f}}}$ for any $t\in{\R}$, 
the winding number of $f(S_r)$ is equal to the difference between 
the signed total changes in angle of the curve $\Phi(t)+i\Psi(t)$ and the curve ${(t+i)}^{\deg{f}}$. 
By this observation, we can calculate the number of roots of $f$ contained in $B_r$ with the help of Kronecker's theorem. 
\begin{cor}
\textit{Suppose that $f(z)$ has no roots on $S_r$. 
Let $w(t)$ denotes the number of sign changes in the Sturm sequence of $\Phi(t)$ and $\Psi(t)$. 
Then, the number of roots of $f$ contained in $B_r$ equals to $\dfrac{w(M)-w(-M)+\deg{f}}{2}$.} 
\end{cor}

For any real polynomial $f$, the sign of $f(t)$ for sufficiently large (resp. small) $t\in{\R}$ 
is determined by the leading coefficient (resp. multiplied by $(-1)^{\deg{f}}$). 
Therefore, in order to determine $w(M)$, we only see the leading coefficients of the Sturm sequence of $\Phi(t)$ and $\Psi(t)$. 
For the rest of the paper, $w(\infty)$ (resp. $w(-\infty)$) denotes the number of sign changes in the leading coefficients (resp. multiplied by $(-1)^{\deg{f_i}}$) of the Sturm sequence.

\subsection{Method for deciding the distribution of roots of a real polynomial.}
Suppose $f(z)$ be a real polynomial of one complex variable $z$. 
Then, we can determine the distribution of roots of $f$ as follows. 

If we want to know the number of real roots of $f$ contained in the closed interval $[a, b]$, then 

1. Check that $a$ and $b$ are not roots of $f$. 

2. Calculate the Sturm sequence of $f(t)$ and $f'(t)$. 

3. By using Sturm's theorem, $w(a)-w(b)$ is equal to the number of real roots of $f$ contained in $[a, b]$. 

If we want to know the number of roots of $f$ contained in $B_r$, then 

1. Calculate the two real polynomials $\Phi(t)$ and $\Psi(t)$ by substituting $z(t)=r\dfrac{t-i}{t+i}$ for $f(z)$. 

2. Check that $f(z)$ has no roots on $S_r$. 
For example, if the resultant of $\Phi(t)$ and $\Psi(t)$ does not equal to $0$, then $f(z)$ has no roots on $S_r$. 

3. Calculate the Sturm sequence of $\Phi(t)$ and $\Psi(t)$. 

4. By Corollary 1 and the definition of $w(\infty)$ and $w(-\infty)$, the number of roots of $f$ contained in $B_r$ is equal to $\dfrac{w(\infty)-w(-\infty)+\deg{f}}{2}$.

\section{Construction of infinite series of non-simple ideal hyperbolic Coxeter polytopes}
In this Section, 
we construct infinite series $\{P_n\}_{n\in{\N}}$ of non-simple ideal hyperbolic Coxeter 4-polytopes 
by glueing ideal hyperbolic Coxeter 4-pyramids along their isometrical facets. 
First, we introduce the vertical projection $p_{\infty}$ from $\infty$ to $\R^3$ and describe how to see hyperbolic 4-polytopes in terms of the projection. 
Second, we review hyperbolic Coxeter 4-pyramids $P_1$ over the product of three simplexes which are completely classified by Tumarkin \cite{T} 
and then construct the infinite family $\{P_n\}_{n\in{\N}}$. 
Finally, we determine the combinatorial structure of $P_n$ in order to calculate the growth rate $\tau_{P_n}$. 
In the sequel, we call 2-faces of 4-polytope {\it faces}.

\subsection{The vertical projection from $\infty$}
First of all, we recall horospheres in $\H^4$. 
A \textit{horosphere $\Sigma$ based at a point at infinity $u$} is defined to be 
a 3-dimensional Euclidean sphere in $\H^4$ tangent to $\R^3$ at $u$ (resp. a Euclidean hyperplane parallel to $\R^3$) 
if $u$ is situated on $\R^3$ (resp. $u=\infty$). 
If we restrict the hyperbolic metric on the horosphere $\Sigma$, it makes a model of 3-dimensional Euclidean geometry. 
\begin{lem} (Theorem6.4.5, \cite{Ra})
Suppose that $P=\cap_{i=1}^{m}H^{-}_{i}$ is a non-compact hyperbolic $4$-polytope and $u$ is a vertex at infinity of $P$. 
Let $\Sigma$ be a horosphere based at $u$ such that $\Sigma$ intersects with $P$ only at the bounding hyperplanes incident to $u$. 
Then, $L(u):=P\cap{\Sigma}$ has the following properties.
\begin{itemize}
\item $L(u)$ is a $3$-dimensional Euclidean polytope in $\Sigma$. 
\item For any bounding hyperplane $H_i$ incident to $u$, $H_i\cap{L(u)}$ is a bounding hyperplane of $L(u)$ in $\Sigma$. 
\item If 2 facets $F_i:=H_i\cap{P}$ and $F_j:=H_j\cap{P}$ make the face of $P$, then the intersection of $F_i\cap{L(u)}$ and $F_j\cap{L(u)}$ is an edge of $L(u)$ and the dihedral angle $\angle{F_i\cap{F_j}}$ is equal to the dihedral angle $\angle{\bigl(F_i\cap{L(u)}\bigr)\cap{\bigl(F_j\cap{L(u)}\bigr)}}$
\end{itemize}
\end{lem}

\begin{figure}[htbp]
\begin{center}
 \includegraphics[scale=.4, clip]{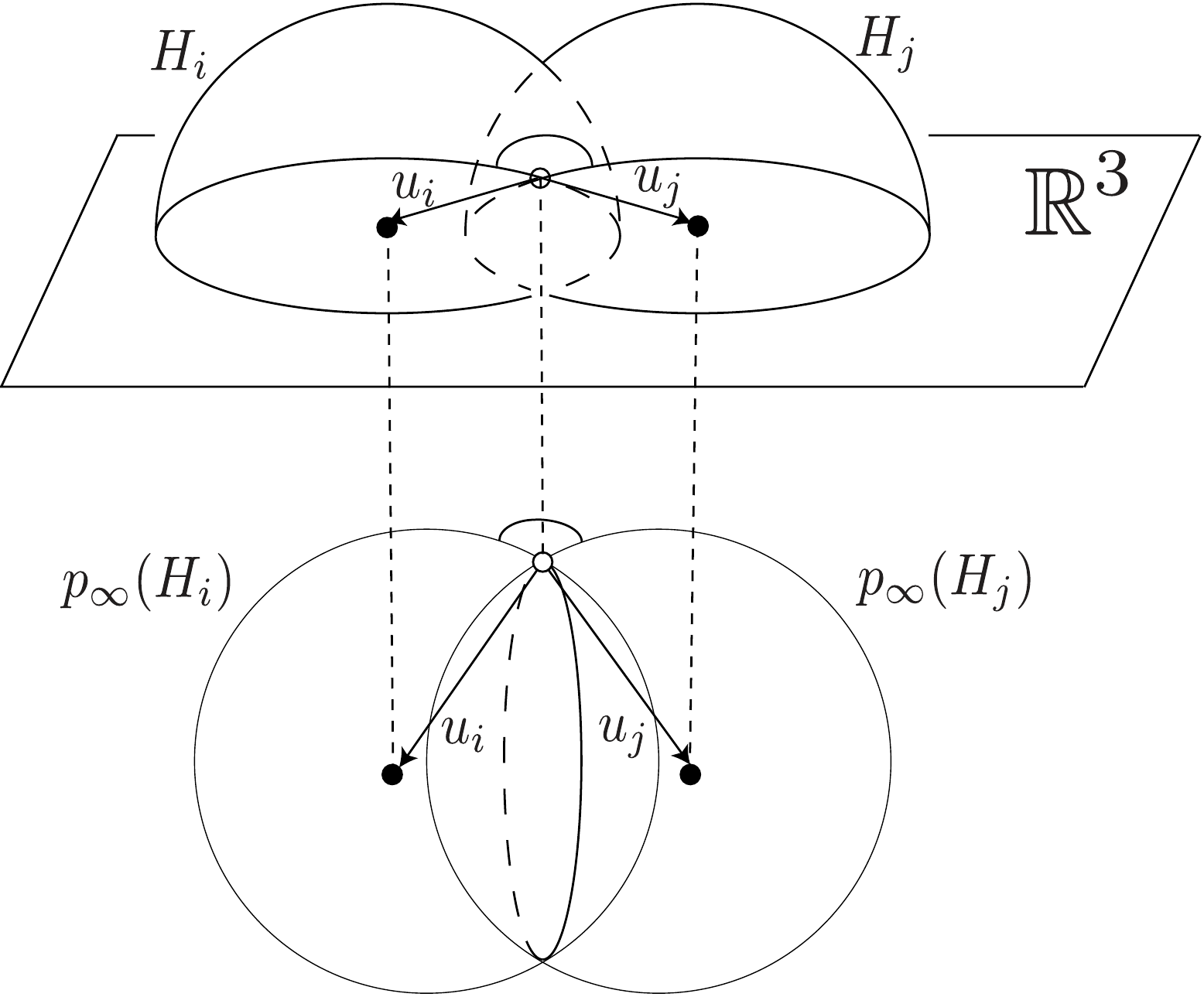}
\end{center}
\caption{The dihedral angle in $\R^3$}
\label{dihedral angle}
\end{figure}

We call the following mapping $p_{\infty}$ \textit{the vertical projection from $\infty$}. 
\[
p_{\infty}: \H^4\to{\R^3}; (x, y, z, t)\mapsto{(x, y, z)}. 
\]
Let $P=\cap_{i=1}^{m}{H^{-}_i}$ be a non-compact hyperbolic 4-polytope and $u$ be a vertex at infinity of $P$. 
By using the translation on $\R^3$ which maps $u$ to $0$ and the inversion with respect to the unit sphere in $\R^4$, we may assume that $u$ is $\infty$. 
If a hyperplane $H_i$ is incident to (resp. not incident to) $\infty$, then $H_i$ is a Euclidean hyperplane (resp. hemisphere) in $\H^4$ orthogonal to $\R^3$. 
Note that any closed half-space $H^{-}_i$ contains $\infty$. 
Since the vertical projection $p_{\infty}$ maps any horosphere $\Sigma$ based at $\infty$ isometrically onto $\R^3$, by using Lemma 2, 
we can treat dihedral angles between 2 bounding hyperplanes of $P$ incident to $\infty$ as corresponding dihedral angles in the $3$-dimensional Euclidean polytope $p_{\infty}(L(\infty))$. 
Suppose that bounding hyperplanes $H_i$ and $H_j$ are not incident to $\infty$. 
By choosing a point in $H_i\cap{H_j}\cap{\R^3}$ and considering the outer normal vectors $u_i$ and $u_j$, we can see the dihedral angle $\angle{H_i\cap{H_j}}$ in $\R^3$ (see Fig \ref{dihedral angle}).

\subsection{Projective image of the ideal hyperbolic Coxeter pyramid $P_1$. }
\begin{theo}( Lemma 10, 11 \cite{T})
There exists an ideal hyperbolic Coxeter 4-pyramid $P_1$ such that the Coxeter scheme $X_1$ of $P_1$ is represented as in Figure \ref{the Coxeter diagram}. 
\begin{figure}[htbp]
\begin{center}
 \includegraphics [scale=.4, clip]{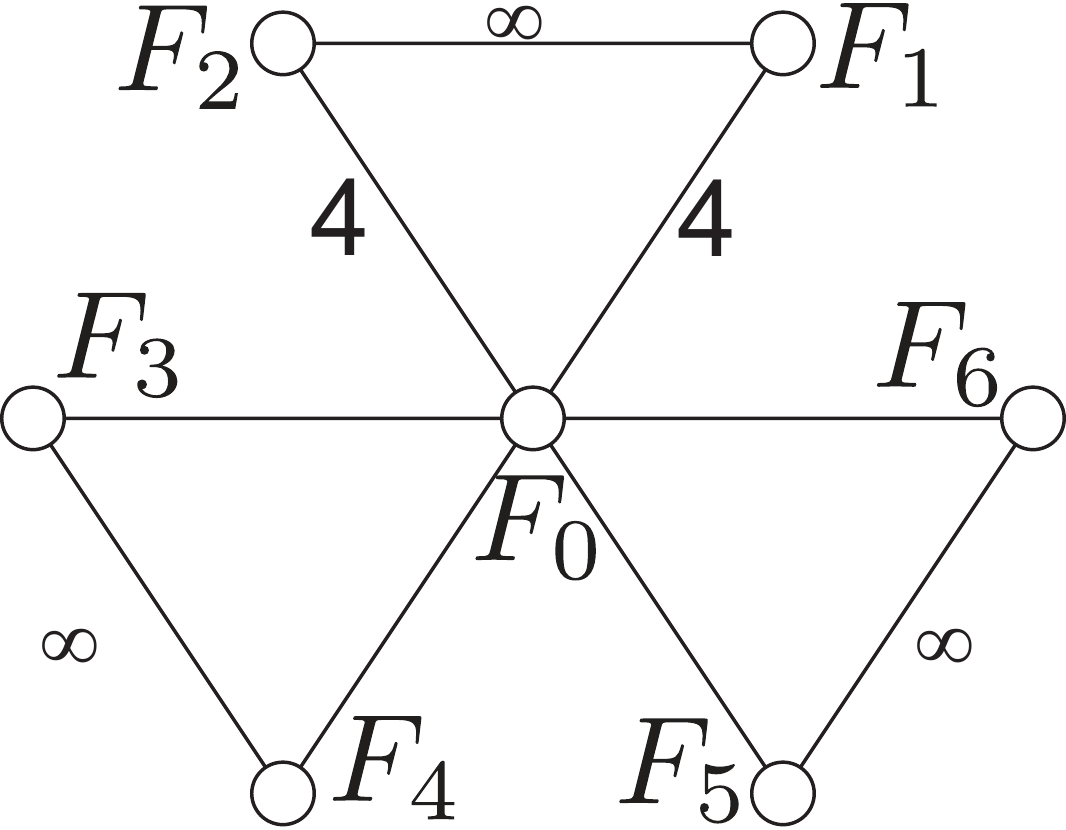}
\end{center}
\caption{The Coxeter scheme $X_1$} 
\label{the Coxeter diagram}
\end{figure} 
\end{theo}

In this Subsection, we use the following notation. 
\begin{itemize}
\item $F_0$ denotes the cubical facet of $P_1$. 
\item The pyramidal facets of $P_1$ are denoted by $F_1, \cdots, F_6$ with the following property : 
$F_i$ and $F_{i+1} (i=1, 3, 5)$ meet at the non-simple vertex of $P_1$ and the dihedral angle of $F_i$ and $F_0$ is equal to $\frac{\pi}{4}$ for $i=1, 2$.   
\item If the intersection of facets $F_i$ and $F_j$ is a face of $P_1$, we denote the face by $f_{ij}$.
\item The non-simple vertex of $P_1$ is denoted by $u$. 
\item The bounding hyperplane of $F_i$ is denoted by $H_i$. 
\end{itemize}

\begin{figure}[htbp]
\begin{center}
 \includegraphics [scale=.4, clip]{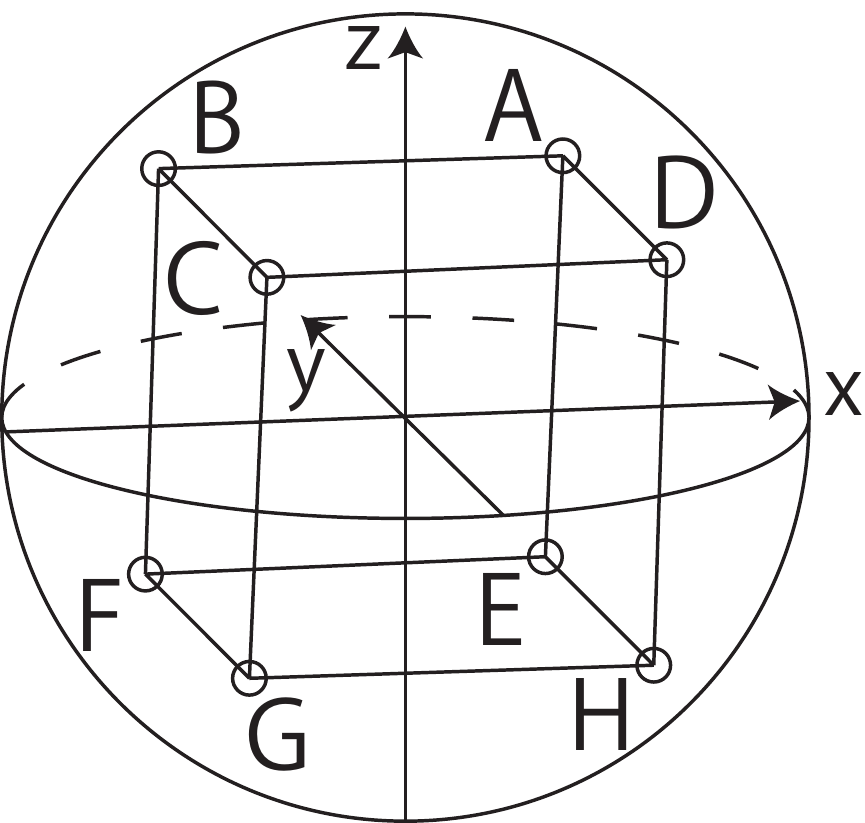}
\end{center}
\caption{$p_{\infty}(P_1)$}
\label{projection1}
\end{figure}

Since the vertex link of $u$ is a Euclidean right rectangular prism, by using isometries of $\H^4$, $P_1$ can be normalized as follows: 
\begin{itemize}
\item The vertex $u$ is $\infty$. 
\item The bounding hyperplane $H_0$ is the unit hemisphere centered at origin. 
\item The bounding hyperplanes $H_1$ and $H_2$ are orthogonal to the $x$-axis. 
\item The bounding hyperplanes $H_3$ and $H_4$ are orthogonal to the $y$-axis. 
\item The bounding hyperplanes $H_5$ and $H_6$ are orthogonal to the $z$-axis. 
\end{itemize}
Under the normalization of $P_1$, we can see $p_{\infty}(P_1)$ as Figure \ref{projection1}, where the coordinates of the  eight points $A, B, C, D, E, F, G$ and $H$ are \\
\hspace{5mm} $A=(\dfrac{1}{\sqrt{2}}, \dfrac{1}{2}, \dfrac{1}{2})$ \hspace{5mm} $B=(-\dfrac{1}{\sqrt{2}}, \dfrac{1}{2}, \dfrac{1}{2})$ \hspace{5mm} $C=(-\dfrac{1}{\sqrt{2}}, -\dfrac{1}{2}, \dfrac{1}{2})$ \hspace{5mm} $D=(\dfrac{1}{\sqrt{2}}, -\dfrac{1}{2}, \dfrac{1}{2})$ \\
\hspace{3mm} $E=(\dfrac{1}{\sqrt{2}}, \dfrac{1}{2}, -\dfrac{1}{2})$ \hspace{2mm} $F=(-\dfrac{1}{\sqrt{2}}, \dfrac{1}{2}, -\dfrac{1}{2})$ \hspace{2mm} $G=(-\dfrac{1}{\sqrt{2}}, -\dfrac{1}{2}, -\dfrac{1}{2})$ \hspace{2mm} $H=(\dfrac{1}{\sqrt{2}}, -\dfrac{1}{2}, -\dfrac{1}{2})$. 

In Figure \ref{projection1}, bounding hyperplanes for quadrangular faces ADHE, ABFE and ABCD are $p_{\infty}(H_1), p_{\infty}(H_3)$ and $p_{\infty}(H_5)$. 
We take a copy of $P_1$, denoted by $P'_1$, and then glue two isometric 4-pyramids $P_1$ and $P'_1$ along the facet $F_1$ of $P_1$ and the facet $F_2$ of $P'_1$. 
\begin{figure}[htbp]
\begin{center}
 \includegraphics [scale=.4, clip]{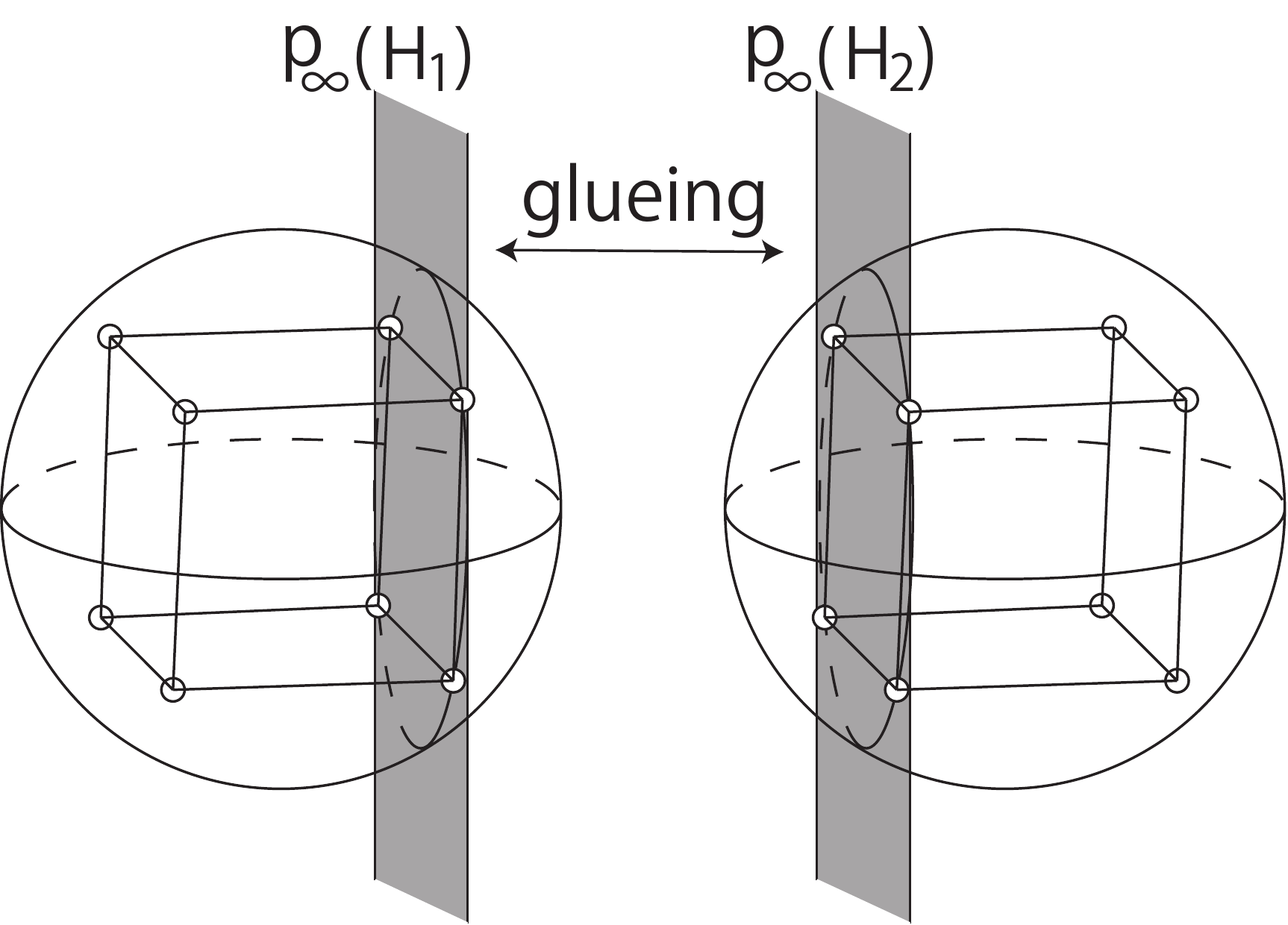}
\end{center}
\caption{The projective images of $P_1$ and $P'_1$}
\end{figure}

\begin{figure}[htbp]
\begin{center}
 \includegraphics [scale=.4, clip]{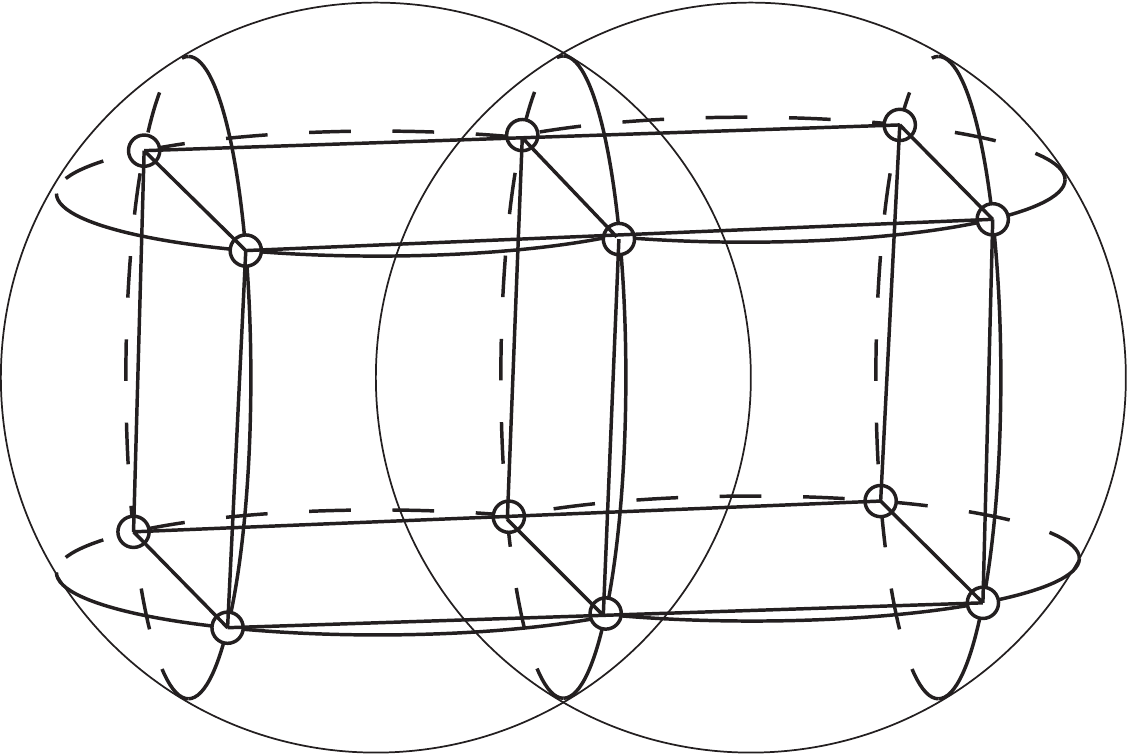}
\end{center}
\caption{the projective image of the resulting 4-polytope.}
\label{projection2}
\end{figure}

Then, we can see the projective image of the resulting 4-polytope $P_2$ as in the Figure \ref{projection2}. 
By the glueing procedure, facets $F_1$ of $P_1$ and $F_2$ of $P'_1$ disappear from $P_2$. 
Since hyperplanes $p_{\infty}(H_3), p_{\infty}(H_4), p_{\infty}(H_5)$ and $p_{\infty}(H_6)$ of $P_1$ and $P'_1$ coincide with each other, 
faces $f_{13}, f_{14}, f_{15}, f_{16}$ in $P_1$ and $f_{23}, f_{24}, f_{25}, f_{26}$ in $P'_1$ also disappear from $P_2$. 
On the other hand, $P_2$ has some new faces; one is the quadrangular face composed by each cubical facet $F_0$ in $P_1$ and $P'_1$ and the other new faces are composed by unions of $f_{34}, f_{45}, f_{56}$ and $f_{63}$ of $P_1$ and $P'_1$. 
Since the facets $F_2$ in $P_1$ and $F_1$ in $P'_1$ do not contribute to the glueing procedure, $P_2$ has the two pyramidal facets $F_1$ and $F_2$.

By summarizing this observation, we see the combinatorial data of $P_2$ as follows. 
\begin{itemize}
\item $P_2$ has 8 facets; 2 cubical facets, 2 pyramidal facets and 4 facets with 6 faces. 
\item $P_2$ has 23 faces; (i) 8 triangular faces come from $F_2$ of $P_1$ and $F_1$ of $P'_1$, (ii) 10 quadrangular faces come from $F_0$ in $P_1$ and $P'_1$, (iii) only one quadrangular face comes from the intersection of $F_1$ in $P_1$ and $F_2$ in $P'_1$, (iv) 4 quadrangular faces come from the union of  $f_{34}, f_{45}, f_{56}$ and $f_{63}$ of $P_1$ and $P'_1$.
\item $P_2$ has 28 edges. 
\item $P_2$ has 13 vertices; only one vertex $\infty$ is non-simple. 
\end{itemize}
Since the two pyramidal facets of $P_2$ are isometric to pyramidal facets $F_1$ and $F_2$ of $P_1$, 
we can repeat this procedure by glueing $P_1$ and $P_2$ along isometric pyramidal facets, and the resulting 4-polytope is denoted by $P_3$. 
This observation implies that we can continue this procedure over and over again. 
The ideal hyperbolic 4-polytope obtained by glueing $n$ copies of $P_1$ along isometric facets $F_1$ and $F_2$ is denoted by $P_n$.

\subsection{Combinatorial structure of $P_n$. }
\begin{lem}
$P_n$ has the following combinatorial data. 
\begin{itemize}
\item[(Facet)] $(n+6)$ facets; $n$ cubical facets, 2 pyramidal facets and the other 4 facets have $(n+4)$ faces.
\item[(Faces)] $(5n+13)$ faces; 8 triangular faces, $5n+1$ quadrilateral faces and 4 $(n+2)$-gonal faces. 
\item[(Edges)] $(8n+12)$ edges. 
\item[(Vertices)] $(4n+5)$ vertices; $4n+4$ simple vertices and only one non-simple vertex. 
\end{itemize}
\end{lem}
\proof By considering the projective image of $P_n$, we can see the assertion. 
Indeed, the projective image $P_n$ consists of $n$ right quadrangular prisms inscribed in closed balls of radius 1 (see Fig 8). 
\begin{figure}[htbp]
\begin{center}
 \includegraphics [scale=.4, clip]{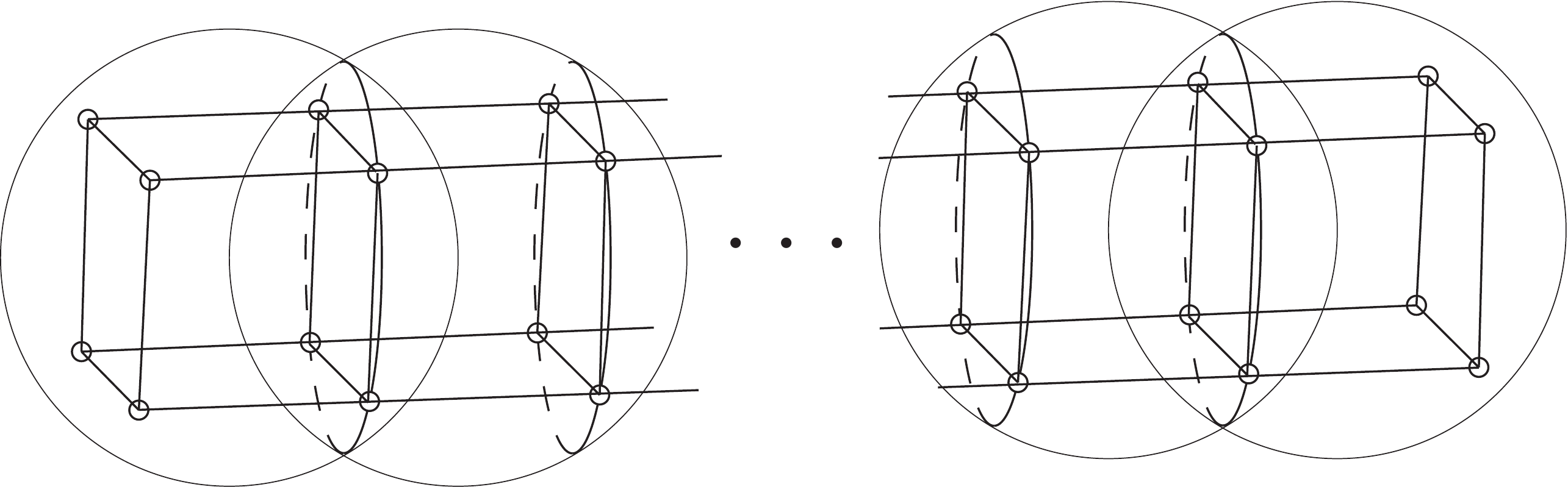}
\end{center}
\caption{The projective image of $P_n$}
\label{projective image of $P_n$}
\end{figure}
\qed.

We use the following notation and terminology in this section. 
\begin{itemize}
\item 2 pyramidal facets of $P_n$ are denoted by $F_1, F_2$.
\item $n$ cubical facets of $P_n$ are denoted by $C_1, \cdots, C_n$. 
Moreover, we suppose that $C_1\cap{F_1}, C_n\cap{F_2}$ and $C_i\cap{C_{i+1}}$ are quadrilateral faces.  
\item The other facets of $P_n$ are denoted by $G_1, G_2, G_3, G_4$. Moreover, we suppose that $G_i\cap{G_{i+1}}$ ($i {\rm \ mod \ } 4$) is a $(n+2)$-gonal face.  
\item $X_n$ denotes the Coxeter scheme of $P_n$. 
\item If a face of $P_n$ has dihedral angle $\frac{\pi}{m}$, we call it a face with $\frac{\pi}{m}$. 
\end{itemize}
\begin{figure}[htbp]
\begin{center}
 \includegraphics [scale=.4, clip]{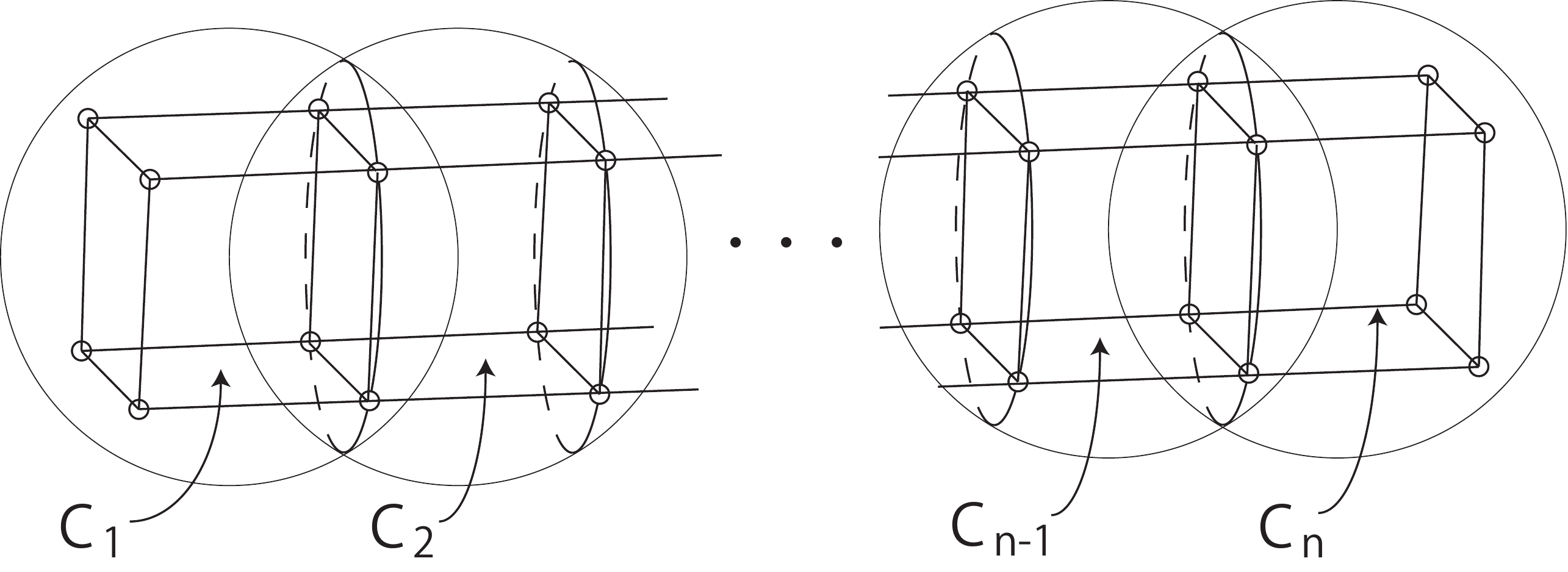}
\end{center}
\caption{The notations for $P_n$}
\label{projective image of $P_n$}
\end{figure}
Let us now determine the elliptic or parabolic subschemes of $X_n$. 

(1) By Lemma 3, $X_n$ has n+6 vertices. 

(2) Since each quadrilateral face $C_i\cap{C_{i+1}}$ is an intersection of glueing facets, 
its dihedral angle $\angle{C_i\cap{C_{i+1}}}$ is equal to $\frac{\pi}{2}$. 
If we glue $P_{n-1}$ and $P_1$ along their isometric pyramidal facets, 
then every faces of $P_{n-1}$ and $P_1$ which is not adjacent to glueing facets does not change. 
Therefore, 
\begin{itemize}
\item the triangular faces $F_i\cap{G_j}$ are faces with $\frac{\pi}{2}$.
\item the $(n+2)$-gonal faces $G_i\cap{G_{i+1}}$ are faces with $\frac{\pi}{2}$. 
\item the quadrilateral faces $G_i\cap{C_j}$ are faces with $\frac{\pi}{3}$. 
\item the quadrilateral faces $C_1\cap{F_1}$ and $C_n\cap{F_2}$ are faces with $\frac{\pi}{4}$. 
\end{itemize}

(3) Each edge of $P_n$ is expressed as the intersection of the three facets. 
\begin{itemize}
\item If an edge is expressed as the intersection of $F_i\cap{G_j}\cap{G_{j+1}}$, 
it corresponds to the elliptic subscheme $A_1\times{A_1}\times{A_1}$ of $X_n$. 
\item If an edge is expressed as the intersection of $F_1\cap{G_i}\cap{C_1}$ or $F_2\cap{G_i}\cap{C_n}$, 
it corresponds to the elliptic subscheme $B_3$ of $X_n$. 
\item If an edge is expressed as the intersection of $G_i\cap{G_{i+1}}\cap{C_j}$, 
it corresponds to the elliptic subscheme $A_3$ of $X_n$. 
\item If an edge is expressed as the intersection of $G_i\cap{C_j}\cap{C_{j+1}}$, 
it corresponds to the elliptic subscheme $A_3$ of $X_n$. 
\end{itemize} 

(4) Each vertex corresponds to the parabolic subscheme of $X_n$. 
\begin{itemize}
\item If a vertex is the intersection of $F_1\cap{G_i}\cap{G_{i+1}}\cap{C_1}$ or $F_2\cap{G_i}\cap{G_{i+1}}\cap{C_n}$, 
it corresponds to the parabolic subscheme $\tilde{B}_3$ of $X_n$. 
\item If a vertex is the intersection of $G_{i}\cap{G_{i+1}}\cap{C_j}\cap{C_{j+1}}$, 
it corresponds to the parabolic subscheme $\tilde{A}_3$ of $X_n$.
\item If a vertex is non-simple, 
it corresponds to the parabolic subscheme $\tilde{A}_1\times{\tilde{A}_1}\times{\tilde{A}_1}$ of $X_n$.
\end{itemize}

\section{The growth function of $P_n$} 
By combining with the combinatorial data of $P_n$ and Steinberg's formula, the growth function $f_n(t)$ of $P_n$ can be calculated as follows.
 
\[
\dfrac{1}{f_n(t^{-1})}=1-\dfrac{n+6}{[2]}+\dfrac{n+11}{[2,2]}+\dfrac{4n}{[2,3]}+\dfrac{2}{[2,4]}-\dfrac{8}{[2,2,2]}-\dfrac{8}{[2,4,6]}-\dfrac{8n-4}{[2,3,4]}. 
\]

By using Mathematica, the growth function $f_n(t)$ can be expressed as 
\[
f_n(t^{-1}):=\dfrac{P_n(t)}{D_n(t)}
\]
where 
\begin{eqnarray*}
P_n(t) &=& (t+1)^3(t^2+1)(t^2-t+1)(t^2+t+1) \\
D_n(t) &=& t^9-(n+3)t^8-(n-4)t^7+(2n-8)t^6 \\ 
&& \ \ \ \ \ \ \ +(2n+8)t^5+(2n-8)t^4-(2n-11)t^3+(3n-5)t^2+(3n+4)t-4(n+1). 
\end{eqnarray*}

\begin{lem}
\textit{All the roots of $D_n(t)$ are simple. }
\end{lem}
\proof We show that the resultant $R(D_n(t), D'_n(t))$ of $D_n(t)$ and $D'_n(t)$ is not equal to $0$ for any $n\in{\N}$. 
By using Mathematica, we can calculate it as follows: 
\begin{eqnarray*}
R(D_n(t), D'_n(t))&=&9367548196608 n^{16}-84315693201408 n^{15}-3211145218356480 n^{14}\\
 && -13452086684085248
   n^{13}-76883986729280512 n^{12}-221310749589989376 n^{11} \\
   &&-369276695931527424
   n^{10}-436823682353681408 n^9-375744535536699392 n^8\\
   &&-227155659791212544n^7-100271146222672128 n^6-28147372028425216 n^5\\
   &&-2791806794781440 n^4-1194005028478976
   n^3-23952968404992 n^2\\
   &&-2787725279232 n
\end{eqnarray*}
By using the Descartes rule \cite{P} (Corollary 1, p.28), $R(D_n(t), D'_n(t))$ has at most one real positive roots as a real polynomial of a real variable $n$. 
We can check the following equalities by using Mathematica. 
\begin{eqnarray*}
R(D_{25}(t), D'_{25}(t)) &=& -5236764089528548306162419869100800. \\
R(D_{26}(t), D'_{26}(t)) &=& 18356309345841539117459400503775232.
\end{eqnarray*}
Hence, $R(D_n(t), D'_n(t))\neq{0}$ for any $n\in{\N}$.\qed

\subsection{The distribution of real roots of $D_n(t)$}
\begin{lem} 
Let $w(t)$ be the number of sign changes in the Sturm sequence of $D_n(t)$ and $D'_n(t)$. 
Then, 
$w(0)=\begin{cases}
        6 & (1\leq{n}\leq{25}) \\
        5 & (26\leq{n})
        \end{cases}
, w(\infty)=\begin{cases} 
                  3 & (1\leq{n}\leq{25}) \\
                  2 & (26\leq{n})
                     \end{cases}$. \\
Moreover, by using Sturm's theorem, the number of real positive roots of $D_n(t)$ is equal to $3$ for any $n\in{\N}$. 
\end{lem}
\proof The equality $D_n(0)=-4(n+1)$ implies that $0$ is not a root of $D_n(t)$ for any $n\in{\N}$. 
By using Mathematica, the Sturm sequence of $D_n(t)$ and $D'_n(t)$ can be calculated and listed in Appendix. 
Let us denote the Sturm sequence of $D_n(t)$ and $D'_n(t)$ as $d_0, \cdots, d_9$ and the $i$-th coefficient of $d_k(t)$ as $a_i^{(k)}$, that is, 
\[
d_k(t)=\sum_{i=0}^{9-k}a_i^{(k)}t^i. 
\]
Then, $w(0)$ (resp. $w(\infty)$) is equal to the number of sign changes in the sequence $a_0^{(0)}, \cdots, a_0^{(9)}$ (resp. $a_9^{(0)}, a_8^{(1)}, \cdots, a_1^{(8)}, a_0^{(9)}$). 
The sign of each coefficient $a_i^{(k)}$ depends on $n\in{\N}$. 
From now on, we determine its signs. 
For example, we consider the sign of $a_0^{(5)}$. 
The sign of $a_0^{(5)}$ depends on the following polynomial $p(n)$ (see Appendix); 
\begin{eqnarray*}
p(n)&=&13008n^8+20600n^7-1607896n^6+2420092n^5+2017855n^4\\ 
 && \hspace{2cm}+899112n^3+1122697n^2-1476508n-45088.
\end{eqnarray*}
Let us first calculate the difference between $p(n+1)$ and $p(n)$. 
\begin{eqnarray*}
p(n+1)-p(n) &=& 52032n^7+254212n^6-4243164n^5-5193210n^4 \\
 && \hspace{2cm}+781934n^3+7841885n^2+7857749n+1704480
\end{eqnarray*}
By the Descartes rule \cite{P} (Corollary 1, p.28), the number of positive real zeroes of $p(n+1)-p(n)$ is at most 2. 
\begin{eqnarray*}
p(2)-p(1) &=& 9055918>0. \\
p(3)-p(2) &=& -140899954<0. \\
p(8)-p(7) &=& -10316213144<0. \\
p(9)-p(8) &=& 16414574600>0.
\end{eqnarray*}
This observation shows that 
\[
 \begin{cases}
   p(2)>p(1) \\
   p(2)>p(3)>\cdots>p(7)>p(8) \\
   p(8)<p(9)<\cdots<p(n)<p(n+1)<\cdots. 
\end{cases}
\]
Moreover, 
\begin{eqnarray*}
p(1) &=& 3363872. \\
p(3) &=& -260324200. \\
p(9) &=& -39144733360. \\
p(10) &=& 162088321532. \\
\end{eqnarray*}
Therefore, we can determine the sign of $a_0^{(5)}$ as follows; 
\[
 a_0^{(5)} \begin{cases}
                   >0 & (n=1, 2) \\
                   <0 & (3\leq{n}\leq{9}) \\
                   >0 & (n\geq{10}). 
                 \end{cases}
 \]
The case of other $a_i^{(k)}$ is considered by analogy, so that we obtain the assertion. \qed

We can calculate $w(-\infty)$ analogously to the proof of Lemma 5. 
\begin{eqnarray*}
w(-\infty) &=& \begin{cases}
                     6 & (1\leq{n}\leq{25}) \\
                     7 & (26\leq{n}). 
                     \end{cases}
\end{eqnarray*}
Therefore, by combining Lemma 5 and Sturm's theorem, we obtain the following proposition.  
\begin{prop}
\textit{The denominator polynomial $D_n(t)$ has the real roots as follows; 
\[ \begin{cases}
\textrm {three positive roots and no negative roots} & (1\leq{n}\leq{25}) \\
\textrm {three positive roots and two negative roots}& (n\geq{26})
\end{cases}
\]}
\end{prop}

\subsection{The distribution of complex roots of $D_n(t)$}
By applying the method in section 3.3, we can verify an upper bound of the absolute values of all complex roots of $D_n(t)$. 

1. Calculate the two real polynomials $\Phi(t)$ and $\Psi(t)$ which satisfy the following identity
\[
D_n(z(t))=\dfrac{\Phi(t)+i\Psi(t)}{(t+i)^{\deg{D_n}}}
\]
where $z(t)=2\dfrac{t-i}{t+i}$. 
By using Mathematica, $\Phi(t)$ and $\Psi(t)$ can be written as follows: 
\begin{eqnarray*}
\Phi(t) &=& -(162n+56)t^9+(6456n-6512)t^7-(2476n-49792)t^5 \\
 && \hspace{2cm} -(7176n+60048)t^3+(894n+13752)t. \\
\Psi(t) &=& (2034n-456)t^8-(8280n-24880)t^6-(7188n+67136)t^4 \\
 && \hspace{2cm} +(4136n+36816)t^2-(14n+2808). 
\end{eqnarray*}

2. By using Mathematica, we can show that the resultant of $\Phi(t)$ and $\Psi(t)$ is not equal to $0$ for any $n\in{N}$.
Therefore $D_n(t)$ has no roots on the circle $S_2$ of radius $2$ centered at the origin. 

3. By using Mathematica, the Sturm sequence of $\Phi(t)$ and $\Psi(t)$ can be calculated. 

4. In a manner similar to the argument in section 5.1, we can calculate the numbers of sign changes $w(\infty)$ and $w(-\infty)$ in the Sturm sequence $\Phi(t)$ and $\Psi(t)$. 
\begin{lem}
\textit{$w(\infty)=8, w(-\infty)=1$ for any $n\in{\N}$. 
By Corollary 1, the number of roots of $Q_n(t)$ contained in the closed disk of radius 2 centered at the origin in the complex plane $\C$ is equal to 8. }
\end{lem}

\begin{theo}
\textit{The growth rate of $P_n$ is a Perron number for any $n\in{\N}$. }
\end{theo}
\proof By Lemma 6, the absolute values of 8 roots of $D_n(t)$ are strictly less than 2. 
Since $\deg{D_n(t)}=9$, if we prove that $D_n(t)$ has a positive real root which is greater than 2, the assertion follows. 
In order to prove that, we consider the number of real roots of $D_n(t)$ which is greater than 2. 
This number is calculated by applying the method in section 3.3, we can see that 
\[
w(2)=\begin{cases} 
       4 & (1\leq{n}\leq{25}) \\ 
       3 & (26\leq{n}).
     \end{cases}
\]
Therefore, by Sturm's theorem, the number of real roots of $D_n(t)$ which is greater than 2 is equal to 1 for any $n\in{\N}$. \qed

\section{Appendix: the Sturm sequence of $D_n(t)$ and $D'_n(t)$} 
In this Section, we list the Sturm sequence considered in Section 5.1. 
\begin{eqnarray*}
d_0(t) &=& t^9-(n+3)t^8-(n-4)t^7+(2n-8)t^6 \\ 
&& \ \ \ \ \ \ \ +(2n+8)t^5+(2n-8)t^4-(2n-11)t^3+(3n-5)t^2+(3n+4)t-4(n+1) \\
d_1(t) &=& 9t^8-8(n+3)t^7-7(n-4)t^6+6(2n+8)t^5 \\ 
&& \ \ \ \ \ \ \ +5(2n+8)t^4+4(2n-8)t^3-3(2n-11)t^2+2(3n-5)t+(3n+4) \\
d_2(t) &=& \dfrac{1}{81} \Bigl\{(8 n^2+66 n)t^7+(7 n^2-61 n+132)t^6+(-12 n^2-60 n-144)t^5 \\
&& \ \ \ \ \ \ \ +(-10 n^2-160 n+240)t^4+(-8 n^2+116 n-498)t^3+(6 n^2-204 n+216)t^2 \\
&& \ \ \ \ \ \ \ \ \ \ +(-6 n^2-224 n-258)t-3 n^2+311 n+312\Bigr\} \\
d_3(t) &=& \dfrac{81}{4n^2(4n+33)^2}\Bigl\{(39 n^4+266 n^3-2673 n^2-1848 n-1936)t^6\\
&& \ \ \ \ \ \ \ \ \ \ \ \ \ \ \ \ \ \ \ +(36 n^4+612 n^3+3956 n^2+4480 n+2112)t^5\\
&& \ \ \ \ \ \ \ \ \ \ \ \ \ \ \ \ \ \ \ +(54 n^4+470n^3-1872 n^2-4372 n-3520)t^4\\
&& \ \ \ \ \ \ \ \ \ \ \ \ \ \ \ \ \ \ \ +(-88 n^4-776 n^3+3866 n^2+6246 n+7304)t^3\\
&& \ \ \ \ \ \ \ \ \ \ \ \ \ \ \ \ \ \ \ +(150 n^4+1374 n^3-3216 n^2-1660 n-3168)t^2\\
&& \ \ \ \ \ \ \ \ \ \ \ \ \ \ \ \ \ \ \ +(162 n^4+2508 n^3+8540 n^2+8870 n+3784)t\\
&& \ \ \ \ \ \ \ \ \ \ \ \ \ \ \ \ \ \ \ -259 n^4-3428 n^3-7161 n^2-8548 n-4576\Bigr\} 
\end{eqnarray*}

From now on, we list the coefficients of the polynomials $d_4(t), \cdots, d_8(t)$. 
\begin{eqnarray*}
\text{The denominator of $f_4(t)$} &=& 81(1936+n(1848+n(2673-n(266+39n))))^2
\end{eqnarray*}
\begin{eqnarray*}
a_5^{(4)}&=&8n^2(4n+33)^2(270 n^6-930 n^5-59765 n^4-72316 n^3-51247 n^2-34920 n+11920)\\
a_4^{(4)}&=&-16n^2(4n+33)^2(51n^6+1630n^5+7368n^4-68445n^3-3176n^2-41152n+16768)\\
a_3^{(4)}&=&8 n^2 (4 n+33)^2 (471 n^6+6452 n^5-5086 n^4-176746 n^3-54403 n^2-120344 n-8944)\\
a_2^{(4)}&=&16 n^2 (4 n+33)^2 (153 n^6-411 n^5-32385 n^4-33106 n^3-44007 n^2-20216 n-7664)\\
a_1^{(4)}&=&-8 n^2 (4 n+33)^2 (579 n^6+14834 n^5+101041 n^4+47610 n^3+25760 n^2+3472 n-25280)\\
a_0^{(4)}&=&16 n^2 (33 + 4 n)^2 (10304 + 60992 n + 92088 n^2 + 112317 n^3 +78944 n^4 + 5932 n^5 + 33 n^6)\\
\end{eqnarray*}
\begin{eqnarray*}
\text{The denominator of $f_5(t)$} &=& 4n^2(33+4n)^2(11920-34920n-51247n^2-72316n^3 \\
&& \hspace{3cm}-59765n^4-930n^5+270n^6)^2
\end{eqnarray*}
\begin{eqnarray*}
a_4^{(5)}&=&-81(39 n^4+266 n^3-2673 n^2-1848 n-1936)^2(246 n^8-5794 n^7+360959 n^6\\
&&+5606880 n^5-3313218 n^4+6140122 n^3-3491843 n^2+2584756n-544176)\\
a_3^{(5)}&=&162(39 n^4+266 n^3-2673 n^2-1848 n-1936)^2 (5289 n^8+5992 n^7-788952 n^6\\
&&-810030 n^5-5107313 n^4+118907 n^3-2823408 n^2+1353973n-43828)\\
a_2^{(5)}&=&-81(39 n^4+266 n^3-2673 n^2-1848 n-1936)^2 (8442 n^8-32742 n^7-1868957 n^6\\
&&-1946748 n^5-4253223 n^4-1203496 n^3-1818280 n^2+440564n-127008)\\
a_1^{(5)}&=&-162(39 n^4+266 n^3-2673 n^2-1848 n-1936)^2 (6261 n^8+27352 n^7-543939 n^6\\
&&+1168425 n^5-740209 n^4-333809 n^3-454006 n^2-793981n+269220)\\
a_0^{(5)}&=&81(39 n^4+266 n^3-2673 n^2-1848 n-1936)^2(13008 n^8+20600 n^7-1607896 n^6\\
&&+2420092 n^5+2017855 n^4+899112 n^3+1122697 n^2-1476508n-45088)
\end{eqnarray*}
\begin{eqnarray*}
\text{The denominator of $f_6(t)$} &=& 81(-1936-1848n-2673n^2+the 266n^3+39n^4)^2 \\
&&(-544176+2584756n-3491843n^2+6140122n^3-3313218n^4 \\
&&+5606880n^5+360959n^6-5794n^7+246n^8)^2
\end{eqnarray*}
\begin{eqnarray*}
a_3^{(6)}&=&-8 n^2 (4 n+33)^2(270 n^6-930 n^5-59765 n^4-72316 n^3-51247 n^2-34920 n+11920)^2\\
&& (403481 n^{10}+2480778 n^9-37969219 n^8-158119702n^7-1100390746 n^6-216055166 n^5\\
&&-1160964773 n^4+282443786 n^3-329580155 n^2+172728524 n-35052620)\\
a_2^{(6)}&=&16 n^2 (4 n+33)^2(270 n^6-930 n^5-59765 n^4-72316 n^3-51247 n^2-34920 n+11920)^2\\
&&(169494 n^{10}+14649 n^9-18830064 n^8+62828800n^7-387398843 n^6+226406803 n^5\\
&&-413299018 n^4+245275527 n^3-138927361 n^2+67186063 n-4007124) \\
a_1^{(6)}&=&8 n^2 (4 n+33)^2(270 n^6-930 n^5-59765 n^4-72316 n^3-51247 n^2-34920 n+11920)^2\\
&&(474903 n^{10}+4516538 n^9-11601465 n^8+104831670n^7+294141284 n^6-180768204 n^5\\
&&+111338775 n^4-296355112 n^3+31452859 n^2-39181768 n+10452012)\\
a_0^{(6)}&=&-16 n^2 (4 n+33)^2 (270 n^6-930 n^5-59765 n^4-72316 n^3-51247 n^2-34920 n+11920)^2\\
&&(252601 n^{10}+1535932 n^9-10172760 n^8+137682333n^7+130244020 n^6+208421539 n^5\\
&&+143139607 n^4+2115857 n^3+44003972 n^2-41200307 n+18745192) 
\end{eqnarray*}
\begin{eqnarray*}
\text{The denominator of $f_7(t)$} &=& 4n^2(33+4n)^2(11920-34920n-51247n^2-72316n^3 \\
&&-59765n^4-930n^5+270n^6)^2(-35052620+172728524n\\
&&-329580155n^2+282443786n^3-1160964773n^4-216055166n^5\\
&&-1100390746n^6-158119702n^7-37969219n^8+2480778n^9\\
&&+403481n^{10})^2
\end{eqnarray*}
\begin{eqnarray*}
a_2^{(7)}&=&81(39 n^4+266 n^3-2673 n^2-1848 n-1936)^2(246 n^8-5794 n^7+360959 n^6+5606880 n^5\\
&&-3313218 n^4+6140122 n^3-3491843 n^2+2584756n-544176)^2(48400755 n^{12}\\
&&+245803454 n^{11}-4721345357 n^{10}-11572421870 n^9-124324436353 n^8-146160412422 n^7\\
&&-206861074257n^6-134297550268 n^5-66775078001 n^4-24225751096 n^3+3620403819 n^2\\
&&-813838328 n+111404496)\\
a_1^{(7)}&=&162(39 n^4+266 n^3-2673 n^2-1848 n-1936)^2(246 n^8-5794 n^7+360959 n^6+5606880 n^5\\
&&-3313218 n^4+6140122 n^3-3491843 n^2+2584756n-544176)^2(9127365 n^{12}\\
&&+43738914 n^{11}-1050600669 n^{10}-2134594907 n^9-221052668 n^8+8764159647 n^7\\
&&+11937399782 n^6+16709700491n^5+4028829086 n^4+2954840024 n^3-2598459169 n^2\\
&&-405956928 n-67672272)\\
a_0^{(7)}&=&-81(39 n^4+266 n^3-2673 n^2-1848 n-1936)^2(246 n^8-5794 n^7+360959 n^6+5606880 n^5\\
&&-3313218 n^4+6140122 n^3-3491843 n^2+2584756n-544176)^2(59130903 n^{12}\\
&&+320783028 n^{11}-5921870437 n^{10}-16668405100 n^9-117418503841 n^8-151967821848 n^7\\
&&-180213457131n^6-140644288440 n^5-51131969275 n^4-32331152680 n^3+5676560341 n^2\\
&&-2814520288 n-23940048) \\
\end{eqnarray*}
\begin{eqnarray*}
\text{The denominator of $f_8(t)$} &=& 81(-1936-1848n-2673n^2+the 266n^3+39n^4)^2 \\
&&(-544176+2584756n-3491843n^2+6140122n^3-3313218n^4 \\
&&+5606880n^5+360959n^6-5794n^7+246n^8)^2(111404496-813838328n\\
&&+3620403819n^2-24225751096n^3-66775078001n^4-134297550268n^5\\
&&-206861074257n^6-146160412422n^7-124324436353n^8\\
&&-11572421870n^9-4721345357n^{10}+245803454n^{11}+48400755n^{12})^2
\end{eqnarray*}
\begin{eqnarray*}
a_1^{(8)}&=&16 n^2 (4 n+33)^2(270 n^6-930 n^5-59765 n^4-72316 n^3-51247 n^2-34920 n+11920)^2\\
&&(403481 n^{10}+2480778 n^9-37969219 n^8-158119702
   n^7-1100390746 n^6-216055166 n^5\\
&&-1160964773 n^4+282443786 n^3-329580155 n^2+172728524 n-35052620)^2\\
&&(1462545045 n^{14}-10472627469n^{13}-402243294759 n^{12}-1104112693071 n^{11}\\
&&-8571517376059 n^{10}-16797900884717 n^9-22904507347277 n^8-22168784110521 n^7\\
&&-14235620251809n^6-6907194126551 n^5-2062300172501 n^4-196719185377 n^3\\
&&-72614586920 n^2+4391952 n-226865664)\\
\end{eqnarray*}
\begin{eqnarray*}
a_0^{(8)}&=&-16 n^2 (4 n+33)^2(270 n^6-930 n^5-59765 n^4-72316 n^3-51247 n^2-34920 n+11920)^2\\
&&(403481 n^{10}+2480778 n^9-37969219 n^8-158119702n^7-1100390746 n^6-216055166 n^5\\
&&-1160964773 n^4+282443786 n^3-329580155 n^2+172728524 n-35052620)^2\\
&&(682442280 n^{14}-13967744415n^{13}-318617986273 n^{12}-866028050552 n^{11}\\
&&-5973136686946 n^{10}-11470936502501 n^9-15278417145211 n^8-15018314214172 n^7\\
&&-9591556809634n^6-5038052836203 n^5-1582742665577 n^4-286371055374 n^3\\
&&-76587929392 n^2-3723242592 n-226865664)\\
\end{eqnarray*}
\begin{eqnarray*}
\text{The numerator of }a_0^{(9)}&=&81(39 n^4+266 n^3-2673 n^2-1848 n-1936)^2(246 n^8-5794 n^7\\
&&+360959 n^6+5606880 n^5-3313218 n^4+6140122 n^3-3491843 n^2\\
&&+2584756n-544176)^2(48400755 n^{12}+245803454 n^{11}-4721345357 n^{10}\\
&&-11572421870 n^9-124324436353 n^8-146160412422 n^7-206861074257n^6\\
&&-134297550268 n^5-66775078001 n^4-24225751096 n^3+3620403819 n^2\\
&&-813838328 n+111404496)^2(36591985143 n^{15}-329358176568n^{14}\\
&&-12543536009205 n^{13}-52547213609708 n^{12}-300328073161252 n^{11}\\
&&-864495115585896 n^{10}-1442487093482529 n^9-1706342509194068n^8\\
&&-1467752091940232 n^7-887326796059424 n^6-391684164932313 n^5\\
&&-109950671986036 n^4-10905495292115 n^3-4664082142496 n^2\\
&&-93566282832n-10889551872)\\
\end{eqnarray*}
\begin{eqnarray*}
\text{The denominator of }a_0^{(9)}&=&4 n (4 n+33)^2(270 n^6-930 n^5-59765 n^4-72316 n^3-51247 n^2-34920 n\\
&&+11920)^2(403481 n^{10}+2480778 n^9-37969219 n^8-158119702
   n^7\\
&&-1100390746 n^6-216055166 n^5-1160964773 n^4+282443786 n^3\\
&&-329580155 n^2+172728524 n-35052620)^2(1462545045 n^{14}\\
&&-10472627469
   n^{13}-402243294759 n^{12}-1104112693071 n^{11}\\
&&-8571517376059 n^{10}-16797900884717 n^9-22904507347277 n^8\\
&&-22168784110521 n^7-14235620251809n^6-6907194126551 n^5\\
&&-2062300172501 n^4-196719185377 n^3-72614586920 n^2+4391952 n\\
&&-226865664)^2
\end{eqnarray*}

\section{Acknowledgements}
The author wishes to express his gratitude to Professor Ruth Kellerhals and Professor Yohei Komori for fruitful discussions of ideas of this paper and their helpful comments concerning Sturm's theorem and its applications for the growth rates. 
This work was partially supported by Grant-in-Aid for JSPS Fellows number 17J05206.

\end{document}